\documentclass{article}
 \usepackage{amsmath}

\newtheorem{thm}{Theorem}[section]

\newtheorem{lem}[thm]{Lemma}

\newtheorem{ex}[thm]{Example}

\newcommand{\be}{\begin{equation}}
\newcommand{\ee}{\end{equation}}
\newcommand{\ben}{\begin{enumerate}}
\newcommand{\een}{\end{enumerate}}
\newcommand{\beq}{\begin{eqnarray}}
\newcommand{\eeq}{\end{eqnarray}}
\newcommand{\beqn}{\begin{eqnarray*}}
\newcommand{\eeqn}{\end{eqnarray*}}

\newcommand{\pa}{\partial}

\newcommand{\qed}{\hspace*{\fill}Q.E.D.}  

\begin{document}
\title{On a Class of  Singular Douglas and Projectively flat Finsler
Metrics}
\author{Guojun Yang  }
\date{}
\maketitle
\begin{abstract}
  Singular Finsler metrics, such as Kropina metrics and $m$-Kropina
  metrics,
   have a lot of applications in the real world. In this paper, we
  study a class of  singular Finsler metrics
   defined by a Riemann metric $\alpha$ and 1-form $\beta$ and  characterize those which are
 respectively Douglasian and locally projectively flat in dimension $n\ge 3$ by some equations.
   Our study shows that the main class induced
   is an $m$-Kropina metric plus a linear part on $\beta$.
     For this class with $m\ne -1$,   the local structure of  projectively flat case
    is determined, and
    it is proved that a Douglas $m$-Kropina metric must be Berwaldian and a projectively flat $m$-Kropina metric must be
     locally Minkowskian. It  indicates that the singular case  is
     quite different from  the regular one.


 {\bf MR(2000) subject classification: }
53A20, 53B40
\end{abstract}

\section{Introduction}
       There are
   two important projective invariants in projective Finsler geometry: the Douglas
  curvature ({\bf D}) and the Weyl curvature (${\bf W}^{o}$ in dimension two and
  {\bf W} in higher dimensions) (\cite{Dou}). A Finsler metric is called
  {\it Douglasian} if ${\bf D}=0$. Roughly speaking, a Douglas metric
 is a Finsler metric having the same geodesics as a Riemannian metric.
  A Finsler metric is said to be {\it locally projectively flat} if at every point, there are local
coordinate systems in which geodesics are straight.
  As we know, the locally projectively flat class of Riemannian metrics
  is very limited, nothing but the class
  of constant sectional
  curvature (Beltrami Theorem). However, the class of locally projectively flat Finsler metrics is very
  rich. Douglas metrics form a rich class of Finsler metrics including
locally projectively flat Finsler metrics, and meanwhile there are
many Douglas metrics which are not locally projectively flat.

In this paper, we will concentrate on a special class of Finsler
metrics:  $(\alpha,\beta)$-metrics, and characterize those which
are Douglasian and locally projectively flat under the condition
(\ref{j2}) below. An {\it $(\alpha,\beta)$-metric} is defined by a
Riemannian metric
 $\alpha=\sqrt{a_{ij}(x)y^iy^j}$ and a $1$-form $\beta=b_i(x)y^i$ on a manifold
 $M$, which can be expressed in the following form:
 $$F=\alpha \phi(s),\ \ s=\beta/\alpha,$$
where $\phi(s)$ is a function satisfying certain conditions. It is
known that $F$ is a  regular Finsler metric if $\beta$ satisfies
 $ \|\beta\|_{\alpha} < b_o$ and $\phi(s)$ is  $C^{\infty}$  on $(-b_o,b_o)$ satisfying
 \be\label{j1}
 \phi(s)>0,\quad \phi(s)-s\phi'(s)+(\rho^2-s^2)\phi''(s)>0,
   \quad (|s| \leq \rho <b_o),
 \ee
 where $b_o$ is a positive
constant (\cite{Shen2}). If $\phi(0)$ is not defined or $\phi$
does not satisfy (\ref{j1}),
 then the $(\alpha,\beta)$-metric $F= \alpha \phi
(\beta/\alpha)$ is singular. Singular Finsler metrics have a lot
of applications in the real world (\cite{AIM} \cite{AHM} ). Z.
Shen also introduces singular Finsler metrics in \cite{Shen3}.

  Assume
  $\phi(s)$ is in the following form
   \be\label{j2}
 \phi(s):=cs+s^m\varphi(s),
   \ee
where $c,m$ are constant with $m\ne 0,1$ and $\varphi(s)$ is a
$C^{\infty}$ function on a neighborhood of $s=0$ with
$\varphi(0)=1$, and further for convenience we put $c=0$ if $m$ is
a negative integer. If $m=0$, we have $\phi(0)=1$ and this case
appears in a lot of literatures. When $m\ge 2$ is an integer,
(\ref{j2})
  is equivalent to the following condition
  $$
  \phi(0)=0,\quad
  \phi^{(k)}(0)=0 \ \ (2\le k\le m-1), \quad \phi^{(m)}(0)=m!.
  $$
  Another interesting case is $c=0$ and $\varphi(s)\equiv 1$ in
  (\ref{j2}), and in this case, $F=\alpha\phi(s)$ is called an
  $m$-Kropina metric, and in particular a Kropina metric when $m=-1$.

The case $\phi(0)=1$ has been studied in  a lot of interesting
research papers  (\cite{LSS}--\cite{LS2} \cite{Shen1}
\cite{Shen2}, \cite{Y}--\cite{Y1}). In \cite{LSS} \cite{Shen1},
the authors respectively study and characterize Douglas
$(\alpha,\beta)$-metrics and locally projectively flat
$(\alpha,\beta)$-metrics in dimension $n\geq 3$ and $\phi(0)=1$,
and further,  the present author solves the case $n=2$ and shows
that the two-dimensional case is quite different from the higher
dimensional ones (\cite{Y1}). In singular case, there are some
papers on the studies of $m$-Kropina metrics and Kropina metrics
(\cite{RR} \cite{ShY} \cite{YN} \cite{YO}). Further, in \cite{Y3},
the present author classifies a class of two-dimensional singular
$(\alpha,\beta)$-metrics $F=\alpha\phi(\beta/\alpha)$ with
$\phi(s)$ satisfying the condition (\ref{j2}) which are Douglasian
and locally projectively flat respectively. In this paper we will
solve the singular case under  the condition (\ref{j2}) in higher
dimensions, which shows   that the singular case is quite
different form the regular condition $\phi(0)=1$ (cf. \cite{LSS}
\cite{Shen1}).

\begin{thm}\label{th01}
  Let $F=\alpha \phi(s)$, $s=\beta/\alpha$, be an $n$-dimensional
  $(\alpha,\beta)$-metric on an open subset $U\subset R^n$ ($n\ge 3$), where $\phi$ satisfies (\ref{j2}).
   Suppose $db\ne0$ in $U$ and that $\beta$ is not parallel with respect to
     $\alpha$.
  If $F$ is a Douglas metric, or locally projectively flat, then $F$ must be in the following
  form
  \be\label{j002}
F=c\bar{\beta}+\bar{\beta}^m\bar{\alpha}^{1-m},\ \
(\bar{\alpha}:=\sqrt{\alpha^2+k\beta^2}, \ \bar{\beta}:=\beta),
 \ee
 where $c,k$ are constant. Note that $\bar{\alpha}$ is Riemannian if
 $k>-1/b^2$.
\end{thm}

  If $b=constant$ in Theorem
   \ref{th01}, there are other classes for the metric $F$ (see Theorem \ref{th1} and Theorem \ref{th51}
   below). Theorem \ref{th01} also holds if $n=2$, but there is
   much difference between $n=2$ and $n\ge 3$ when we determine
   the local structures of $F$ in (\ref{j002}) which is Douglasian
   or locally projectively flat (cf. \cite{Y3}).

Theorem \ref{th01}  naturally induces an important class of
singular Finsler metric---$m$-Kropina metric
$F=\beta^m\alpha^{1-m}$. When $m=-1$, $F=\alpha^2/\beta$ is called
a Kropina metric. There have been some research papers on Kropina
metrics (\cite{RR} \cite{YN} \cite{YO}). In \cite{ShY}, the
present author and Z. Shen characterize $m$-Kropina metrics which
are weakly Einsteinian.

\bigskip

Next we   determine the local structure of the
 metric $F=c\beta+\beta^m\alpha^{1-m}$ which are Douglasian and locally projectively flat respectivley when
 $m\ne -1$. The method is the application of the following deformation on
 $\alpha$ and $\beta$ which is defined by
  \be\label{cr71}
\widetilde{\alpha}:=b^m\alpha, \ \
\widetilde{\beta}:=b^{m-1}\beta.
 \ee
The deformation (\ref{cr71}) first appears in  \cite{ShY} for the
research on weakly Einstein $m$-Kropina metrics. It also appears
in \cite{Y3}. It is very useful for $m$-Kropina metrics.
Obviously, if $F$ is an $m$-Kropina metric, then $F$ keeps
formally unchanged, namely,
$$F=\beta^m\alpha^{1-m}=\widetilde{\beta}^m\widetilde{\alpha}^{1-m}.$$
Further, $\widetilde{\beta}$ has unit length with respect to
$\widetilde{\alpha}$, that is,
$||\widetilde{\beta}||_{\widetilde{\alpha}}=1$.

\begin{thm}\label{th001}
 Let $F=c\beta+\beta^m\alpha^{1-m}$ be an $n(\ge 3)$-dimensional Douglas $(\alpha,\beta)$-metric,
 where $c,m$ are constant with $m\ne
 0,\pm 1$. Then  we have the following  cases:

 \ben
\item[{\rm (i)}] {\rm($c=0$)}  $F$ can be written as
 $F=\widetilde{\alpha}^{1-m}\widetilde{\beta}^m$,
  where  $\widetilde{\beta}$
  is parallel with respect to the Riemann metric
  $\widetilde{\alpha}$, and
further $\alpha,\beta$ are related
    with $\widetilde{\alpha},\widetilde{\beta}$ by
    \be\label{ycw17}
  \alpha=\eta^{\frac{m}{m-1}}\widetilde{\alpha}, \ \ \ \beta=\eta\widetilde{\beta},
    \ee
where $\eta=\eta(x)>0$ is a scalar function. Further, $F$ is
actually Berwaldian.

\item[{\rm (ii)}] {\rm($c\ne 0$)} $F$ can be written as
$F=c\eta\widetilde{\beta}+\widetilde{\beta}^m\widetilde{\alpha}^{1-m}$,
where  $\widetilde{\beta}$
  is parallel with respect to the Riemann metric
  $\widetilde{\alpha}$ with
  $\eta\widetilde{\beta}$ being closed. Furhter we have   (\ref{ycw17}).
 \een
\end{thm}

\begin{thm}\label{th04}
 Let $F=c\beta+\beta^m\alpha^{1-m}$ be an $n(\ge 3)$-dimensional
 locally projectively flat $(\alpha,\beta)$-metric, where $c,m$ are constant with $m\ne
 0,\pm 1$. Then  we have the following  cases:
\ben \item[{\rm (i)}] {\rm($c=0$)}  $F$ can be written as
 $F=\widetilde{\alpha}^{1-m}\widetilde{\beta}^m$,
  where $\widetilde{\alpha}$ is flat and $\widetilde{\beta}$
  is parallel with respect to $\widetilde{\alpha}$, and thus $\widetilde{\alpha}$ and
  $\widetilde{\beta}$ can be locally written as
 \be\label{ycw16}
\widetilde{\alpha}=|y|,\ \ \widetilde{\beta}=y^1.
 \ee
Further $\alpha,\beta$ are related
    with $\widetilde{\alpha},\widetilde{\beta}$ by (\ref{ycw17}).
 Moreover $F$ is locally Minkowskian.

\item[{\rm (ii)}] {\rm($c\ne 0$)} $F$ can be written as $
  F=c\eta\widetilde{\beta}+\widetilde{\beta}^m\widetilde{\alpha}^{1-m},
 $
  where (\ref{ycw16}) and (\ref{ycw17}) hold with
  $\eta=\eta(x^1)>0$. In this case, $F$ is Berwaldian, or
 locally Minkowskian if and only if  $c=0$ or $\eta=constant$ in
 (\ref{ycw17}); and here $\eta=constant$ implies $\alpha$ is flat
 and $\beta$ is parallel.
 \een
\end{thm}

For the two-dimensional case, we have proved that the metric
$F=c\beta+\alpha^2/\beta$ is always Douglasian, the $m$-Kropina
metric in Theorem \ref{th001}(i) is locally Minkowskian
(determined by Theorem \ref{th04}(i)), and the metric $F$ in
Theorem \ref{th001}(ii) is locally projectively flat if
additionally $m\ne -3$ (\cite{Y3}). When $\widetilde{\alpha}$ is
Not flat and $\widetilde{\beta}$ is parallel with respect to
$\widetilde{\alpha}$, then the $m$-Kropina metric $F$ in Theorem
\ref{th001}(i) is Douglasian but Not locallly projectively flat,
and a family of concrete examples to this case are given in the
last section.

When $m=-1$,  the deformation (\ref{cr71}) cannot be applied to
Theorem \ref{th001} and \ref{th04} to determine the local
structure of $F=c\beta+\alpha^2/\beta$ which is Douglasian or
locallly projectively flat. See the general characterization in
Theorem \ref{th61} and \ref{th62} respectively below. In
\cite{Y4}, we further prove  that for the dimensions $n\ge 2$, if
$F=c\beta+\alpha^2/\beta$ is locally projectively flat with
constant flag curvature, then $F$ is locally Minkowskian. If
$c\beta$ is small, then
 $F=(\alpha^2+c\beta^2)/\beta=\bar{\alpha}^2/\beta$
 is a  Kropina metric. In \cite{Y5}, the present author has shown some non-trivial examples of
 Kropina metrics which are locally projectively flat.

\bigskip

\noindent{\bf Open Problem:} Determine the local structure of the
$n(\ge 3)$-dimensional metric $F=c\beta+\alpha^2/\beta$ which is
Douglasian or  locallly projectively flat.

\section{Preliminaries}

Let $F=F(x,y)$ be a Finsler metric on an $n$-dimensional manifold
$M$.
      In local coordinates, the spray coefficients $G^i$ are
      defined by
 \beq \label{G1}
 G^i:=\frac{1}{4}g^{il}\big \{[F^2]_{x^ky^l}y^k-[F^2]_{x^l}\big \}.
 \eeq
If $F$ is a Douglas metric, then  $G^i$   are in the following
form:
 \beq\label{G2}
   G^i=\frac{1}{2}\Gamma_{jk}^i(x)y^jy^k+P(x,y)y^i,
 \eeq
where $\Gamma_{jk}^i(x)$ are local functions on $M$ and $P(x,y)$
is a local positively homogeneous function of degree one in $y$.
It is easy to see that $F$
 is a Douglas metric if and only if $G^iy^j-G^jy^i$ is a
 homogeneous polynomial in $(y^i)$ of degree three, which by
 (\ref{G2}) can be written as (\cite{BaMa}),
  $$G^iy^j-G^jy^i=\frac{1}{2}(\Gamma^i_{kl}y^j-\Gamma^j_{kl}y^i)y^ky^l.$$

According to G. Hamel's result, a Finsler metric $F$ is
projectively flat in $U$ if and only if
 \be\label{01}
 F_{x^my^l}y^m-F_{x^l}=0.
 \ee
 The above formula implies that $G^i=Py^i$ with $P$ given by
 \be\label{02}
  P=\frac{F_{x^m}y^m}{2F}.
  \ee

Consider an $(\alpha,\beta)$-metric $F =\alpha \phi
(\beta/\alpha)$.
  The spray coefficients $G^i_{\alpha}$ of
$\alpha$
 are given by
  $$G^i_{\alpha}=\frac{1}{4}a^{il}\big \{[\alpha^2]_{x^ky^l}y^k-[\alpha^2]_{x^l}\big
  \}.$$
Let $\nabla \beta = b_{i|j} y^i dx^j$  denote the covariant
derivatives of $\beta$ with respect to $\alpha$ and define
 $$r_{ij}:=\frac{1}{2}(b_{i|j}+b_{j|i}),\ \ s_{ij}:=\frac{1}{2}(b_{i|j}-b_{j|i}),\ \
 r_j:=b^ir_{ij},\ \ s_j:=b^is_{ij},\ \ s^i:=a^{ik}s_k,$$
 where $b^i:=a^{ij}b_j$ and $(a^{ij})$ is the inverse of
 $(a_{ij})$.
 By (\ref{G1}) again, the spray coefficients $G^i$ of $F$
are given by:
  \be\label{y20}
  G^i=G^i_{\alpha}+\alpha Q s^i_0+\alpha^{-1}\Theta (-2\alpha Q
  s_0+r_{00})y^i+\Psi (-2\alpha Q s_0+r_{00})b^i,
  \ee
where $s^i_j=a^{ik}s_{kj}, s^i_0=s^i_ky^k,
s_i=b^ks_{ki},s_0=s_iy^i$, and
 $$
  Q:=\frac{\phi'}{\phi-s\phi'},\ \
  \Theta:=\frac{Q-sQ'}{2\Delta},\ \
  \Psi:=\frac{Q'}{2\Delta},\ \ \Delta:=1+sQ+(b^2-s^2)Q'.
 $$

By (\ref{y20}) one  can see that $F=\alpha\phi(\beta/\alpha)$ is a
Douglas metric if and only if
 \be\label{y21}
 \alpha Q (s^i_0y^j-s^j_0y^i)+\Psi (-2\alpha
 Qs_0+r_{00})(b^iy^j-b^jy^i)=\frac{1}{2}(G^i_{kl}y^j-G^j_{kl}y^i)y^ky^l,
\ee
  where $G^i_{kl}:=\Gamma^i_{kl}-\gamma^i_{kl}$,  $\Gamma^i_{kl} $ are given in (\ref{G2}) and $
  \gamma^i_{kl}:=\pa^2G^i_{\alpha}/\pa y^k\pa y^l.$

Further,  $F=\alpha\phi(\beta/\alpha)$ is projectively flat on
$U\subset R^n$ if and only if \be
(a_{ml}\alpha^2-y_my_l)G^m_{\alpha}+\alpha^3Qs_{l0}+\Psi\alpha(-2\alpha
 Qs_0+r_{00})(\alpha b_l-sy_l)=0,\label{y21*}
\ee
 where $y_l=a_{ml}y^m$.

\bigskip

The following lemma is obvious.

\begin{lem}\label{lem0}
 If $Q=ks$, where $k$ is a constant, then
 $\phi(s)=c\sqrt{1+ks^2}$ for some constant $c$.
\end{lem}

\section{Equations in a Special Coordinate System}\label{s3}

Fix an arbitrary point $x\in M$ and take  an orthogonal basis
  $\{e_i\}$ at $x$ such that
   $$\alpha=\sqrt{\sum_{i=1}^n(y^i)^2},\ \ \beta=by^1.$$
Then we change coordinates $(y^i)$ to $(s, y^a)$ such that
  $$\alpha=\frac{b}{\sqrt{b^2-s^2}}\bar{\alpha},\ \
  \beta=\frac{bs}{\sqrt{b^2-s^2}}\bar{\alpha}, $$
where $\bar{\alpha}=\sqrt{\sum_{a=2}^n(y^a)^2}$. Let
 $$\bar{r}_{10}:=r_{1a}y^a, \ \ \bar{r}_{00}:=r_{ab}y^ay^b, \ \
 \bar{s}_0:=s_ay^a.$$
We have $\bar{s}_0=b\bar{s}_{10},s_1=bs_{11}=0$. In the following,
we also put
 $$\bar{G}^0_{10}:=G^a_{1b}y_ay^b,\ \bar{G}^1_{00}:=G^1_{ab}y^ay^b,\ \bar{G}^a_{00}:=G^a_{bc}y^by^c,\ etc.$$
Then by the above coordinate $(s,y^a)$ and using (\ref{y21}) and
(\ref{y21*}), it follows from \cite{LSS} \cite{Shen1} we have the
following lemmas:

 \begin{lem}\label{lem1} {\rm (\cite{LSS})}
 For $n\ge 2$, an $(\alpha,\beta)$-metric
 $F=\alpha\phi(\beta/\alpha)$ is a Douglas metric if and only if
 there hold the following four identities:
 \be \label{cr002}
 \frac{bQ\bar{s}^a_0s-\Psi
 r_{11}s^2by^a}{b^2-s^2}\bar{\alpha}^2-\Psi\bar{r}_{00}by^a=
  \frac{s^2}{2(b^2-s^2)}(\bar{G}^a_{10}+\bar{G}^a_{01}-
    G^1_{11}y^a)\bar{\alpha}^2-\frac{1}{2}\bar{G}^1_{00}y^a,
 \ee
 \be\label{cr003}
 \frac{bQs^2s^a_1}{b^2-s^2}\bar{\alpha}^2+(-2\Psi
 s\bar{r}_{10}+2\Psi Qb^2s^1_0-Qs^1_0)by^a =
 \frac{G^a_{11}s^3}{2(b^2-s^2)}\bar{\alpha}^2+
 \frac{1}{2}\big\{\bar{G}^a_{00}-(\bar{G}^1_{10}+\bar{G}^1_{01})y^a\big\}s,
 \ee
 \be\label{cr004}
 \frac{bs}{b^2-s^2}(s^a_1y^b-s^b_1y^a)Q\bar{\alpha}^2=
 \frac{s^2}{2(b^2-s^2)}(G^a_{11}y^b-G^b_{11}y^a)\bar{\alpha}^2+
 \frac{1}{2}(\bar{G}^a_{00}y^b-\bar{G}^b_{00}y^a),
 \ee
 \be\label{cr005}
 (\bar{s}^a_0y^b-\bar{s}^b_0y^a)bQ=\frac{s}{2}\big\{(\bar{G}^a_{10}+
 \bar{G}^a_{01})y^b-(\bar{G}^b_{10}+\bar{G}^b_{01})y^a\big\}.
 \ee
\end{lem}

 \begin{lem}\label{lem2}
 {\rm(\cite{LSS} \cite{Shen1})} $(n\ge 2)$ Let $F=\alpha\phi(\beta/\alpha)$ be an
 $(\alpha,\beta)$-metric.
  Suppose $\Psi$ is dependent on $s$, $\beta$ is not parallel with
 respect to $\alpha$ and $\beta$ is closed.  Then $F$ is a Douglas metric if and only if
 \beq
  b_{i|j}&=&2\bar{\tau} \big\{\delta b_ib_j+\eta
  (b^2a_{ij}-b_ib_j)\big\}.\label{c1}\\
  2\Psi&=&\frac{\lambda s^2+\mu (b^2-s^2)}{\delta s^2+\eta
  (b^2-s^2)},\label{c2}
 \eeq
 where
 $\bar{\tau}=\bar{\tau}(x),\lambda=\lambda(x),\mu=\mu(x),\delta=\delta(x),\eta=\eta(x)$ are
 scalar functions satisfying $\lambda\eta-\mu\delta\ne 0$. $F$ is
 projectively flat if and only if (\ref{c1}), (\ref{c2}) and
  \be\label{c3}
  G^i_{\alpha}=\rho y^i-\bar{\tau} \big\{\lambda \beta^2+\mu
  (b^2\alpha^2-\beta^2)\big\}b^i
  \ee
 hold, where $\rho:=\rho_i(x)y^i$ is a 1-form.
\end{lem}

\begin{lem} \label{lem3}
{\rm(\cite{Shen1})}
For $n\ge 2$, if $s_{ab}=0$, then an $(\alpha,\beta)$-metric
 $F=\alpha\phi(\beta/\alpha)$ is locally projectively flat if and only if
  \beq
 0&=&\bar{G}^a_{10}\bar{\alpha}^2-\bar{G}^0_{10}y^a,\label{cr006}\\
 0&=&(\bar{r}_{00}+\frac{s^2r_{11}\bar{\alpha}^2}{b^2-s^2})b\Psi
 -\frac{s^2}{2(b^2-s^2)}(2\bar{G}^0_{10}-G^1_{11}\bar{\alpha}^2)+\frac{1}{2}\bar{G}^1_{00},\label{cr007}\\
 0&=&\bar{G}^a_{00}-\Big\{\frac{2bQ(1-2b^2\Psi)\bar{s}_{10}}{s}+4b\Psi\bar{r}_{10}
 +2\bar{G}^1_{10}\Big\}y^a+\frac{s(G^a_{11}s-2bQs_{1a})\bar{\alpha}^2}{b^2-s^2}.\label{cr008}
  \eeq
  where $G^i_{jk}$ are the spray coefficients of $\alpha$.
\end{lem}

Note that in (\ref{cr004}), if $\lim_{s\rightarrow 0}sQ=0$ and
$Q/s$ is dependent on $s$, then we can get $s^a_1=0$. The zero
limit is a key factor to prove $\beta$ is closed using
(\ref{cr004}) and (\ref{cr005}). In singular case, we generally
don't have $\lim_{s\rightarrow 0}sQ=0$.

\section{Douglas $(\alpha,\beta)$-metrics}

In this section, we characterize a class of $n(\ge3)$-dimensional singular
$(\alpha,\beta)$-metrics which are Douglas metrics. We have the
following theorem.

 \begin{thm}\label{th1}
  Let $F=\alpha \phi(s)$, $s=\beta/\alpha$, be an $n(\ge3)$-dimensional
  $(\alpha,\beta)$-metric on an open subset $U\subset R^n$, where $\phi$ satisfies (\ref{j2}).
   Suppose that $\beta$ is not parallel with respect to
     $\alpha$.
  Then $F$ is a Douglas metric if and only if one of the following
  cases holds:
  \ben
 \item[{\rm (i)}]  $\phi$ and $\beta$ satisfy
 \be\label{ygjcw}
\phi(s)=cs+\frac{1}{s},\ \ \  s_{ij}=\frac{b_is_j-b_js_i}{b^2},
 \ee
 where $c$ is a constant.

\item[{\rm (ii)}]
  $\phi$ and $\beta$ satisfy
   \beq
    \phi(s)&=&k_1s+s^m(1+k_2 s^2)^{\frac{1-m}{2}},\label{y5} \\
     b_{i|j}&=&2\tau \big\{mb^2a_{ij}-(m+1+k_2b^2)b_ib_j\big\},
     \label{y6}
   \eeq
   where $\tau=\tau(x)$ is a scalar function and $k_1,k_2$ are
constant.

  \item[{\rm (iii)}]  $\phi$ and $\beta$ satisfy
   \beq
    \phi(s)&=&s^m(1+k s^2)^{\frac{1-m}{2}},\label{y16}\\
   r_{ij}&=&2\tau
   \big\{mb^2a_{ij}-(m+1+kb^2)b_ib_j\big\}-\frac{m+1+2kb^2}{(m-1)b^2}(b_is_j+b_js_i),\label{y17}\\
    s_{ij}&=&\frac{b_is_j-b_js_i}{b^2},\label{y0017}
   \eeq
  where $k$ is constant and $\tau=\tau(x)$ is a scalar.
  \een
 \end{thm}

In Theorem \ref{th1} (iii), if $b=constant$, then $k=-1/b^2$ in
(\ref{y16})--(\ref{y17}), and we get
   \beq
    \phi(s)&=&s^m\big\{1-(\frac{s}{b})^2\big\}^{\frac{1-m}{2}},\label{y18}\\
    r_{ij}&=&2\bar{\tau}(b^2a_{ij}-b_ib_j)-\frac{1}{b^2}(b_is_j+b_js_i),\label{y19}
   \eeq
where $\bar{\tau}:=m\tau$. Note that if $n=2$,  (\ref{y19}) is
equivalent to $b=constant$ (see \cite{LS2}), and clearly
(\ref{y0017}) holds automatically.

\subsection{$d\beta=0$}

Assume $\phi$ satisfies (\ref{j2}),  $\beta$ is not parallel with
 respect to $\alpha$ and $\beta$  is closed. Obviously $F$ is not of Randers type. So by Lemma \ref{lem2}
 we have (\ref{c2}). We first determine
$\lambda,\eta,\delta,\mu$ in (\ref{c2}). Rewrite (\ref{c2}) as
follows
 \be\label{y56}
 [\delta s^2+\eta (b^2-s^2)]\phi''=[\lambda s^2+\mu
 (b^2-s^2)][\phi-s\phi'+(b^2-s^2)\phi''].
 \ee
Plug
 $$\phi(s)=a_1s+s^m(1+a_{m+1}s+a_{m+2}s^2+a_{m+3}s^3+a_{m+4}s^4)+o(s^{m+4})$$
into (\ref{y56}). Let $p_i$ be the coefficients of $s^i$ in
(\ref{y56}). First $p_{m-2}=0$ gives
 \be\label{y57}
 \eta=\mu b^2.
 \ee
  Plugging (\ref{y57}) into  $p_m=0$
yields
 \be\label{y58}
 \delta=\lambda b^2-\frac{m+1}{m}\mu b^2.
 \ee

 \noindent{\bf Case A. } Assume $m=-1$. Plug (\ref{y57}),
(\ref{y58}) and $m=-1$ into (\ref{y56}) and then we get
 $$s^2\phi''+s\phi'-\phi=0,$$
 whose solution is given by (\ref{ygjcw}).

\noindent{\bf Case B. } Assume $m\ne -1$.
 Plugging (\ref{y57}) and (\ref{y58}) into  $p_{m+2}=0$ yields
 \be\label{y59}
  \lambda=[m(m-1)+2a_{m+2}b^2]\epsilon, \ \ \mu=m(m-1)\epsilon,
 \ee
where $\epsilon=\epsilon(x)\ne 0$ is a scalar. It is easy to see
that
 \be\label{y059}
 \lambda \eta-\mu\delta=m(m+1)(m-1)^2b^2\epsilon^2\ne 0.
 \ee
 Plug (\ref{y57}), (\ref{y58}) and (\ref{y59}) into
(\ref{c2}) and we get
 \be\label{y60}
 2\Psi=\frac{\phi''}{\phi-s\phi'+(b^2-s^2)\phi''}=\frac{m(m-1)+2a_{m+2}s^2}{m(m-1)b^2+(1-m^2+2a_{m+2}b^2)s^2},
 \ee
which can be rewritten as
 \be\label{y61}
 \phi''=\frac{-m+k_2s^2}{(1+k_2s^2)s^2}(\phi-s\phi'),
 \ee
 where we put
 $$k_1=a_1,\ \ k_2=-2a_{m+2}/(m-1).$$
Solving the differential equation (\ref{y61}) gives (\ref{y5}).
Plug (\ref{y57}), (\ref{y58}) and (\ref{y59}) into (\ref{c1}) and
we get (\ref{y6}), where we put
 \be\label{y060}
 \tau:=(m-1)\epsilon b^2\bar{\tau}.
 \ee

\subsection{$d\beta\ne 0$}
We will deal with the equations (\ref{cr002})--(\ref{cr005})
respectively.

\bigskip

\noindent {\bf Step (1):} By Lemma \ref{lem0} and the assumption
on $\phi$, we see $Q/s$ is dependent on $s$. So by (\ref{cr005}),
we have
 $$\bar{s}^a_0y^b-\bar{s}^b_0y^a=0,$$
 from which we have $s_{ab}=0$ since $n\ge 3$. Therefore, we
 obtain  (\ref{y0017}).

\bigskip

\noindent {\bf Step (2):} We rewrite (\ref{cr004}) as
 \be\label{cr45}
 0=s\big[2b\phi's_{1a}+s(\phi-s\phi')G^a_{11}\big]y^b-s\big[2b\phi's_{1b}+s(\phi-s\phi')G^b_{11}\big]y^a+(b^2-s^2)(\phi-s\phi')\theta_{ab},
 \ee
where $\theta_{ab}$ are defined by
 $$\theta_{ab}:=\bar{G}^a_{00}y^b-\bar{G}^b_{00}y^a.$$
Plug
$$
 \phi(s)=a_1s+s^m(1+a_{m+1}s+a_{m+2}s^2+a_{m+3}s^3+o(s^{m+3})
$$
 into (\ref{cr45}). Let $p_i$ denote the coefficient of $s^i$ in
 (\ref{cr45}). By $p_m=0$ we get
  \be\label{cr46}
 \theta_{ab}=\frac{2m(s_{1a}y^b-s_{1b}y^a)}{(m-1)b}.
  \ee
Substituting (\ref{cr46}) into $p_{m+2}=0$ yields
 \be\label{cr47}
 T_by^a-T_ay^b=0,
 \ee
where $T_a$ are defined by
 $$T_a:=(m-1)^2bG^a_{11}+2(m-m^2+2a_{m+2}b^2)s_{1a}.$$
Since $n\ge 3$, by (\ref{cr47}) we have $T_a=0$, which are written
as
 \be\label{cr48}
 G^a_{11}=-\frac{2(m-m^2+2a_{m+2}b^2)}{(m-1)^2b}s_{1a}.
 \ee
Finally, plug (\ref{cr46}) and (\ref{cr48}) into (\ref{cr004}) and
then we obtain
 \be\label{cr49}
 Q=-\frac{m+ks^2}{(m-1)s},
 \ee
where we have used the fact that $s_{1a}y^b-s_{1b}y^a\ne 0$ since
$\beta$ is not closed and $s_{ab}=0$, and $k$ is defined by
 $$k:=-\frac{2a_{m+2}}{m-1},$$
Solving the ODE (\ref{cr49}) we get $\phi(s)$ given by
(\ref{y16}).

\bigskip

\noindent {\bf Step (3):} Plug (\ref{cr48}) and (\ref{cr49}) into
(\ref{cr003}) and then we get
 $$
 A_0s^2+mb^2A_1=0,
 $$
where $A_0,A_1$ are polynomials in $(y^a)$ independent of $s$.  By
$A_1=0$ we get
 \be\label{cr50}
 \bar{G}^a_{00}=\Big\{\bar{G}^1_{10}+\bar{G}^1_{01}-\frac{2\bar{r}_{10}}{b}
 -\frac{2(m+1)\bar{s}_{10}}{(m-1)b}\Big\}y^a+\frac{2ms_{1a}\bar{\alpha}^2}{(m-1)b}.
 \ee
Plugging (\ref{cr50}) into $A_0=0$ gives
 $$2(m+1)\big[(m-1)\bar{r}_{10}+(m+1+2kb^2)\bar{s}_{10}\big]y^a=0.$$
So if $m\ne -1$ we get
 \be\label{cr51}
\bar{r}_{10}=-\frac{m+1+2kb^2}{m-1}\bar{s}_{10}.
 \ee
 Now we show $m-kb^2\ne 0$ if $m\ne -1$, which will be needed in the
following. If $m-kb^2= 0$, then $b=constant$, and by (\ref{cr51})
get
 $$0=\bar{r}_{10}+\bar{s}_{10}=-\frac{(m+1)\bar{s}_{10}}{m-1},$$
which is impossible since $\bar{s}_{10}\ne 0$.

\bigskip

\noindent {\bf Step (4):} By $s_{ab}=0$ and a simple analysis on
(\ref{cr002}), we see (\ref{cr002}) can be written as
 \be\label{cr52}
 \frac{s^2}{2(b^2-s^2)}(G^1_{11}-\gamma)\delta_{ab}+\frac{1}{4}(G^1_{ab}+G^1_{ba})=
   b\Psi (\frac{r_{11}s^2}{b^2-s^2}\delta_{ab}+r_{ab}),
 \ee
where $\gamma:=G^a_{1a}+G^a_{a1}$ (not summed) which is
independent of the index $a$. By (\ref{cr49}) and the definition
of $\Psi$ we have
 \be\label{cr53}
 \Psi=\frac{ks^2-m}{2\big[(1+m+kb^2)s^2-mb^2\big]}.
 \ee
Plug $s_{ab}=0$ and (\ref{cr53}) into (\ref{cr52}) and we obtain
 $$
B_0s^4+bB_1s^2+mb^3B_2=0,
 $$
where $B_0,B_1,B_2$ are scalar functions independent of $s$. Then
by $B_2=0$ we have
 \be\label{cr54}
G^1_{ba}=\frac{2r_{ab}}{b}-G^1_{ab}.
 \ee
If $m\ne -1$, using $m-kb^2\ne 0$, plug (\ref{cr54}) into
$B_0=0,B_1=0$ and then we obtain
 \be\label{cr55}
 r_{11}=\frac{b(1+kb^2)(\gamma-G^1_{11})}{m-kb^2},\ \
 r_{ab}=\frac{mb(\gamma-G^1_{11})}{m-kb^2}\delta_{ab}.
 \ee

\bigskip

Now summed up from the  above, it follows from (\ref{cr51}) and
(\ref{cr55}) that (\ref{y17}) holds if $m\ne -1$, where $\tau$ is
defined by
 $$
 \tau:=\frac{G^1_{11}-\gamma}{2b(m-kb^2)}.
$$

\subsection{The inverse of the case $m=-1$}

We have shown that if $F=c\beta+\alpha^2/\beta$ is a Douglas
metric, then $s_{ij}$ are given by (\ref{ygjcw}). There are
different ways to show the inverse is also true. We prove the
inverse by (\ref{y21}). We only need to show the left hand side of
(\ref{y21}) are polynomials in $y$ of degree three.

Plug
 $$
 \phi(s)=cs+\frac{1}{s},\ \ s^i_{0}=\frac{b^is_0-\beta s^i}{b^2}
 $$
into the left hand side of (\ref{y21}), and then we get
 $$
 \frac{1}{2b^2}\big\{(y^js^i-y^is^j)(\alpha^2-c\beta^2)+(b^iy^j-b^jy^i)r_{00}\big\},
 $$
which are clearly polynomials in $y$ of degree three.

\section{Projectively flat $(\alpha,\beta)$-metrics}
In this section, we characterize a class of $n(\ge3)$-dimensional singular
$(\alpha,\beta)$-metrics which are projectively flat. We have the following
theorem.

\begin{thm}\label{th51}
  Let $F=\alpha \phi(s)$, $s=\beta/\alpha$, be an $n(\ge3)$-dimensional
  $(\alpha,\beta)$-metric on an open subset $U\subset R^n$, where $\phi$ satisfies (\ref{j2}).
   Suppose that $\beta$ is not parallel with respect to
     $\alpha$.
  Let $G^i_{\alpha}$  be the spray coefficients of $\alpha$.
  Then $F$ is  projectively flat in $U$ with
 $G^i=P(x,y)y^i$ if and only if one of the following cases holds:

 \ben
 \item[{\rm (i)}]  $\phi(s)$ and $\beta$ satisfy
  (\ref{ygjcw}), and  $G^i_{\alpha}$
  satisfy
  \be\label{w001}
G^i_{\alpha}=\rho
y^i-\frac{r_{00}}{2b^2}b^i-\frac{\alpha^2-c\beta^2}{2b^2}s^i.
 \ee
In this case, the projective factor $P$ is given by
 \be\label{w0001}
 P=\rho-\frac{1}{b^2(\alpha^2+c\beta^2)}\Big\{(\alpha^2-c\beta^2)s_0+r_{00}\beta\Big\}.
 \ee

  \item[{\rm (ii)}]  $\phi(s)$ and $\beta$ satisfy (\ref{y5}) and
  (\ref{y6}), and  $G^i_{\alpha}$
  satisfy
 \be\label{w1}
 G^i_{\alpha} =\rho y^i-\tau(m\alpha^2-k_2\beta^2)b^i.
  \ee
 In this case, the projective factor $P$ is given by
 \be\label{w01}
  P=\rho +\tau \alpha
  \Big\{s(-m+k_2s^2)-s^2(1+k_2s^2)\frac{\phi'}{\phi}\Big\}.
  \ee

\item[{\rm (iii)}] $\phi(s)$ and $\beta$ satisfy
(\ref{y16})--(\ref{y0017}), and  $G^i_{\alpha}$
  satisfy
   \be\label{w3}
   G^i_{\alpha}=\rho
   y^i+\Big\{\frac{2k\beta s_0}{(m-1)b^2}-\tau(m\alpha^2-k\beta^2)\Big\}b^i
   -\frac{m\alpha^2+k\beta^2}{(m-1)b^2}s^i.
   \ee
 In this case, the projective factor $P$ is given by
  \be\label{w03}
  P=\rho-2m\tau\beta-\frac{2m}{(m-1)b^2}s_0.
  \ee
  \een
The above function $\rho= \rho_i(x)y^i$  is a 1-form.
\end{thm}

{\it Proof :} Our proof of Theorem \ref{th51} breaks into two
cases: $m=-1$ and $m\ne -1$. Firstly by (\ref{cr006}) we have
 \be\label{cr62}
 G^a_{1b}=\frac{\xi}{2}\delta_{ab},
 \ee
where $\xi=\xi(x)$ is a scalar function.

\bigskip

\noindent {\bf Step 1. } Assume  $m=-1$.

Plug $\phi(s)=cs+1/s$ into (\ref{cr008}) and we get
 \be\label{cr63}
G^a_{11}=-\frac{1-kb^2}{b}s_{1a},\ \
\bar{G}^a_{00}=(\frac{2\bar{r}_{10}}{b}+2\bar{G}^1_{10})y^a-\frac{\bar{\alpha}^2}{b}s_{1a}.
 \ee
Next substitute $\phi(s)=cs+1/s$ and (\ref{cr62}) into
(\ref{cr007}) and we have
 \be\label{cr64}
 \bar{G}^1_{00}=-\frac{\bar{r}_{00}}{b},\ \
 G^1_{11}=\xi-\frac{r_{11}}{b}.
 \ee
Now by (\ref{cr63}) and (\ref{cr64}) we get (\ref{w001}), where
$\rho$ is defined by
 $$
\rho:=\frac{1}{2}\xi y^1+(\frac{r_{1a}}{b}+G^1_{1a})y^a.
 $$

Finally, we solve the projective factor. Plug $\phi(s)=cs+1/s$,
(\ref{w001}) and $s^i_{0}=(b^is_0-\beta s^i)/b^2$ into
(\ref{y20}), and then we get  $G^i=Py^i$ with $P$ given by
(\ref{w0001}).

\bigskip

\noindent {\bf Step 2. } Assume  $m\ne -1$.

\noindent{\bf Case A:} Assume $d\beta=0$. Plugging (\ref{y59})
into (\ref{c3}) gives (\ref{w1}). Next we show (\ref{w01}). By
(\ref{y6}) we have
 \be\label{j61}
 r_{00}=2\tau \big\{m b^2\alpha^2-(1+m+k_2b^2)\beta^2\big\}.
 \ee
Now plug $s_{i0}=0,s_0=0$ and (\ref{w1}), (\ref{y61}) and
(\ref{j61}) into (\ref{y20}), and then we obtain (\ref{w01}).

 \bigskip

\noindent {\bf Case B:} Assume $d\beta\ne 0$. Plugging
(\ref{cr49}) and (\ref{cr51}) into (\ref{cr008}) gives
 \be\label{cr66}
 \bar{G}^a_{00}=\big(2\bar{G}^1_{10}-\frac{4kb\bar{s}_{10}}{m-1}\big)y^a-\frac{2m\bar{\alpha}^2s_{1a}}{(m-1)b},\
 \ G^a_{11}=-\frac{2(m+kb^2)s_{1a}}{(m-1)b}.
 \ee
Next substituting (\ref{cr53}), (\ref{cr55})  and (\ref{cr62})
into (\ref{cr007}) gives
 \be\label{cr67}
 \bar{G}^1_{00}=-2mb\tau \bar{\alpha}^2, \ \
 G^1_{11}=\xi-2b(m-kb^2)\tau.
 \ee
Now by (\ref{cr66}) and (\ref{cr67}) we get (\ref{w3}), where
$\rho$ is defined by
 $$
\rho:=\frac{1}{2}\xi y^1+(-\frac{2kbs_{1a}}{m-1}+G^1_{1a})y^a.
 $$

Finally, we solve the projective factor. By (\ref{y17}) and
(\ref{y0017}) we have
 \be\label{cr68}
 r_{00}=2\tau \big\{m b^2\alpha^2-(1+m+k_2b^2)\beta^2\big\}-2\frac{m+1+2kb^2}{(m-1)b^2}\beta
 s_0, \ \ \ s^i_{0}=\frac{b^is_0-\beta s^i}{b^2}.
 \ee
Now plug (\ref{cr49}), (\ref{cr53}), (\ref{w3})  and (\ref{cr68})
 into (\ref{y20}), and then we obtain (\ref{w03}).

\section{Proof of Theorem \ref{th001} and Theorem \ref{th04}}

Based on Theorem \ref{th1} and Theorem \ref{th51}, we give a
general characterization for $F=c\beta+\beta^m\alpha^{1-m}$ to be
Douglasian and locally projectively flat respectively.

\begin{thm}\label{th61}
  Let $F=c\beta+\beta^m\alpha^{1-m}$ be an $n$-dimensional
  $(\alpha,\beta)$-metric on an open subset $U\subset R^n$ ($n\ge 3$), where $c,m$ are constant with $m\ne 0, 1$.
   Then for some scalar function $\tau=\tau(x)$, we have the following cases:
    \ben
\item[{\rm (i)}] {\rm($m=-1$)} $F$  is a Douglas metric if and
only if $\beta$ satisfies
 \be\label{cr70}
  s_{ij}=\frac{b_is_j-b_js_i}{b^2}.
  \ee

\item[{\rm (ii)}] {\rm($c\ne 0,m\ne -1$)} $F$  is a Douglas metric
if and only if $\beta$ satisfies
 \be\label{y06}
   b_{i|j}=2\tau \big\{mb^2a_{ij}-(m+1)b_ib_j\big\},
   \ee

  \item[{\rm (iii)}] {\rm($c= 0,m\ne -1$)} $F$
is a Douglas metric if and only if $\beta$ satisfies (\ref{cr70})
and
  \be\label{cr69}
   r_{ij}=2\tau \big\{mb^2a_{ij}-(m+1)b_ib_j\big\}
   -\frac{m+1}{(m-1)b^2}(b_is_j+b_js_i),
   \ee
  \een
\end{thm}

\begin{thm}\label{th62}
  Let $F=c\beta+\beta^m\alpha^{1-m}$ be an $n$-dimensional
  $(\alpha,\beta)$-metric on an open subset $U\subset R^n$ ($n\ge 3$), where $c,m$ are constant with $m\ne 0, 1$.
   Then for some scalar function  $\tau=\tau(x)$ and 1-form $\rho=\rho_i(x)y^i$, we have the following cases:
 \ben
\item[{\rm (i)}] {\rm($m=-1$)} $F$  is projectively flat if and
only if $\beta$ satisfies (\ref{cr70}) and $G^i_{\alpha}$ satisfy
 \be\label{cw1}
G^i_{\alpha}=\rho
y^i-\frac{r_{00}}{2b^2}b^i-\frac{\alpha^2-c\beta^2}{2b^2}s^i.
 \ee

\item[{\rm (ii)}] {\rm($c\ne 0;m\ne -1$)} $F$  is projectively
flat if and only if $\beta$ satisfies (\ref{y06}) and
$G^i_{\alpha}$ satisfy
  \be\label{cw3}
 G^i_{\alpha} =\rho y^i-m\tau \alpha^2b^i.
  \ee

\item[{\rm (iii)}] {\rm($c= 0;m\ne -1$)} $F$ is projectively flat
if and only if $\beta$ satisfies (\ref{cr70}) and (\ref{cr69}),
and $G^i_{\alpha}$ satisfy
   \be\label{cw4}
  G^i_{\alpha}=\rho
   y^i-m\tau \alpha^2b^i+
   \frac{m}{(1-m)b^2}\alpha^2s^i.
   \ee

  \een
\end{thm}

We can use the deformation (\ref{cr71}) to simplify (\ref{cr70}),
(\ref{cr69}) and (\ref{cw4}), which is shown as follows:

\begin{lem}\label{lem61}
  For a pair $(\alpha,\beta)$, suppose $\beta$ satisfies (\ref{cr69}) and
(\ref{cr70}).
   Then under the deformation (\ref{cr71}),
$\widetilde{\beta}$ must be parallel with respect to
$\widetilde{\alpha}$. Further, if the spray coefficients
$G^i_{\alpha}$ of $\alpha$ satisfy (\ref{cw4}),
   then $\widetilde{\alpha}$ is projectively flat.
 \end{lem}

 {\it Proof :}
By (\ref{cr69}) and (\ref{cr70}), a direct computation under
(\ref{cr71}) gives $\widetilde{r}_{ij}=0$ and
$\widetilde{s}_{ij}=0$ respectively. Thus $\widetilde{\beta}$ is
parallel with respect to $\widetilde{\alpha}$. If (\ref{cw4})
holds, then under (\ref{cr71}) we have
 $$
\widetilde{G}^i_{\widetilde{\alpha}}=\Big[\rho-2m\tau
\beta-\frac{2ms_0}{(m-1)b^2}\Big] y^i.
 $$
So $\widetilde{\alpha}$ is projectively flat.

We can also give another simple proof for $m\ne -1$. Define
$F:=\beta^m\alpha^{1-m}$. If (\ref{cr69}) and (\ref{cr70}) hold,
then $F$ is a Douglas metric by Theorem \ref{th61}(iii). Since $F$
keeps formally unchanged under (\ref{cr71}), by Theorem
\ref{th61}(iii) we have
 \beq
   \widetilde{r}_{ij}&=&2\widetilde{\tau} \big\{m\widetilde{a}_{ij}-(m+1)\widetilde{b}_i\widetilde{b}_j\big\}
   -\frac{m+1}{m-1}(\widetilde{b}_i\widetilde{s}_j+\widetilde{b}_j\widetilde{s}_i),\label{cr72}\\
 \widetilde{s}_{ij}&=&\widetilde{b}_i\widetilde{s}_j-\widetilde{b}_j\widetilde{s}_i,\ \ \widetilde{b}^2=1.\label{cr73}
   \eeq
Contracting (\ref{cr72}) by $\widetilde{b}^i$ and then by
$\widetilde{b}^j$ and using $\widetilde{r}_i+\widetilde{s}_i=0$,
it is easy to get $\widetilde{\tau}=0$, $\widetilde{r}_{ij}=0$ and
$\widetilde{s}_i=0$. So by (\ref{cr73}) we have
$\widetilde{s}_{ij}=0$.  Further, if
 $G^i_{\alpha}$ satisfy (\ref{cw4}), then $F$ is locally
 projectively flat by Theorem \ref{th62}(iii). Again, since $F$
keeps formally unchanged under (\ref{cr71}), by Theorem
\ref{th62}(iii) and $\widetilde{b}=1$, and using
$\widetilde{\tau}=0$ and $\widetilde{s}_i=0$, we have
 $$
 \widetilde{G}^i_{\widetilde{\alpha}}=\widetilde{\rho}
   y^i-m\widetilde{\tau} \widetilde{\alpha}^2\widetilde{b}^i+
   \frac{m}{1-m}\widetilde{\alpha}^2\widetilde{s}^i=\widetilde{\rho}
   y^i,
 $$
 which imply $\widetilde{\alpha}$ is projectively flat.     \qed

\bigskip

{\it Proof of Theorem \ref{th001} :}

\bigskip

 If $c=0$, then we have
$F=\beta^m\alpha^{1-m}$. Since $F=\beta^m\alpha^{1-m}$  with $m\ne
-1$ is a Douglas metric, by Theorem \ref{th61}(iii) we have
(\ref{cr69}) and (\ref{cr70}). Put
$\eta:=||\beta||_{\alpha}^{1-m}$ and then we get (\ref{ycw17}).
Then by Lemma \ref{lem61} we complete the proof of Theorem
\ref{th001}(i). If $c\ne 0$, then $\beta$ is closed  by Theorem
\ref{th61}(ii). Thus $\eta\widetilde{\beta}$ is closed. This
completes the proof of Theorem \ref{th001}(ii).  \qed

\bigskip

{\it Proof of Theorem \ref{th04} :}

\bigskip

\noindent {\bf Case A:} Assume $c=0$. Since
$F=\beta^m\alpha^{1-m}$ is a Douglas metric, we have (\ref{cr69}),
 (\ref{cr70}) and (\ref{cw4}) by Theorem \ref{th62}(iii). Then
 Lemma \ref{lem61} shows that under the deformation (\ref{cr71}),
 $\widetilde{\beta}$ is parallel with respect to
 $\widetilde{\alpha}$, and $\widetilde{\alpha}$ is projectively
 flat. Thus we can first locally express $\widetilde{\alpha}$ in the
 following form
  \be\label{cr75}
\widetilde{\alpha} =\frac{\sqrt{(1+\mu |x|^2)|y|^2-\mu\langle
 x,y\rangle^2}}{1+\mu|x|^2},
 \ee
where $\mu$ is the constant sectional curvature of
$\widetilde{\alpha}$. Since $\widetilde{\beta}=\widetilde{b}_iy^i$
is of course a closed 1-form which is conformal with respect to
$\widetilde{\alpha}$, it has been shown in \cite{Yu} the following
  \be\label{cr76}
 \widetilde{b}_i=\frac{kx^i+(1+\mu|x|^2)e_i-\mu\langle e,x\rangle
 x^i}{(1+\mu|x|^2)^{\frac{3}{2}}},\ \ \
 \widetilde{b}^i=\sqrt{1+\mu|x|^2}(kx^i+e_i).
 \ee
where $k$ is a constant and $e=(e_i)$ is a constant vector, and
$\widetilde{b}_i=\widetilde{a}_{ij}\widetilde{b}^j$. By
(\ref{cr76}) we have
 \be\label{cr77}
 1=\widetilde{b}^2=||\widetilde{\beta}||^2_{\widetilde{\alpha}}=|e|^2+\frac{k^2|x|^2+2k \langle e,x\rangle -\mu \langle
 e,x\rangle^2}{1+\mu |x|^2}.
 \ee
 It is easy to conclude from (\ref{cr77}) that $\mu=0$. So
 $\widetilde{\alpha}$ is flat. Thus $\widetilde{\alpha}$ and
 $\widetilde{\beta}$ can be locally expressed as (\ref{ycw16}).

\bigskip

\noindent {\bf Case B:} Assume $c\ne 0$. In this case, we only
need to require additionally that $\beta$ be closed by Theorem
\ref{th62}(ii). Then since $\beta=\eta \widetilde{\beta}=\eta y^1$
is closed, we see $\eta=\eta(x^1)$. Now we can easily verify that
for the metric
  $F=c\eta\widetilde{\beta}+\widetilde{\beta}^m\widetilde{\alpha}^{1-m}$,
  (\ref{01}) holds.
  So $F$ is projectively flat with $G^i=Py^i$. Further, by
  (\ref{02}) we can get the projective factor $P$ given by
 \be\label{ycw107}
 P=\frac{c\eta_1}{2F}(y^1)^2, \ \ \eta_1:=\eta_{x^1},
 \ee
and the scalar flag curvature $K$ is given by
 \be\label{ycw108}
K=
\frac{c(y^1)^3}{2F^3}\Big\{\frac{3c\eta_1^2y^1}{2F}-\eta_{11}\Big\},\
 \ \eta_{11}:=\eta_{x^1x^1}.
\ee
 Then by (\ref{ycw107}) and (\ref{ycw108}), $F$ is Berwaldian, or
 locally Minkowskian if and only if  $c=0$ or $\eta=constant$.

\section{A local representation}
 We have show that if $\widetilde{\alpha}$ is Not flat and
$\widetilde{\beta}$ is parallel with respect to
$\widetilde{\alpha}$, then the $m$-Kropina metric $F$ in Theorem
\ref{th001}(i) is Douglasian but Not locallly projectively flat.
In this section, we give a family of examples to this case.

Firstly we show a lemma based on \cite{Y} (also see \cite{Y2}).

\begin{lem}\label{lem71}
 Let $\widetilde{\alpha}$ be an $n$-dimensional Riemann metric which is
 locally conformally flat, and $\widetilde{\beta}$ is a 1-form. Then $\widetilde{\beta}$
 is a Killing form  $\widetilde{r}_{ij}=0$ with unit length if and only if $\widetilde{\alpha}$ and
 $\widetilde{\beta}$ can be locally expressed as
  \be\label{cr80}
 \widetilde{\alpha}=\sqrt{\frac{|y|^2}{|u|^2}},\ \ \widetilde{\beta}=\frac{\langle
 u,y\rangle}{|u|^2},
  \ee
  where $u:=(u^1(x),\cdots,u^n(x))$ is a vector satisfying the
  following PDEs (fixed $i,j$):

  \be\label{cr81}
  \frac{\pa u^i}{\pa x^j}+\frac{\pa u^j}{\pa
  x^i}=0 \ {\rm (} \forall \thinspace i\ne  j{\rm )},\ \ \ \ \
  \frac{\pa u^i}{\pa x^i}=\frac{\pa u^j}{\pa
  x^j} \ {\rm (} \forall \thinspace  i, j{\rm )}.
 \ee
 Further, if $n=2$, $\widetilde{\beta}$ is parallel with respect to $\widetilde{\alpha}$,
 and $\widetilde{\alpha}$ is flat.

\end{lem}

{\it Proof :} Suppose $\widetilde{\beta}$ satisfies
   $\widetilde{r}_{ij}=0$ and has  unit length. Since $\widetilde{\alpha}$ is locally conformally
   flat, we can express it as $\widetilde{\alpha}=e^{\frac{1}{2}\sigma(x)}|y|$. In this case, firstly we
   can express $\widetilde{\beta}=e^{\sigma}\langle
 u,y\rangle$, and then by \cite{Y}, $u$ satisfies (\ref{cr81}).
 Since $\widetilde{\beta}$ has unit length, clearly we have
 $e^{\sigma}=1/|u|^2$. So we get (\ref{cr80}).

 Conversely, suppose $\widetilde{\alpha}$ and $\widetilde{\beta}$ are given by
 (\ref{cr80}) with $u$ satisfying (\ref{cr81}). Clearly $\widetilde{\beta}$
 has unit length. Next we verify $\widetilde{r}_{ij}=0$. It has been shown in
 \cite{Y} that if $\widetilde{\alpha}$ and $\widetilde{\beta}$
 are given by (\ref{cr80}) with $u$ satisfying (\ref{cr81}), then $\widetilde{\beta}$ is a conformal
 form satisfying
 \be\label{cr82}
 \widetilde{r}_{ij}=\frac{\pa u^1}{\pa x^1}+\frac{1}{2}u^k\sigma_k,
 \ee
where
$$\sigma:=-\ln(|u|^2),\ \ \sigma_k:=\sigma_{x^k}.$$
Then (\ref{cr82}) becomes
 \be\label{cr83}
 \widetilde{r}_{ij}=\frac{1}{|u|^2}\big(\frac{\pa u^1}{\pa x^1}|u|^2-u^iu^k\frac{\pa u^i}{\pa
 x^k}\big).
 \ee
By (\ref{cr81}), we have
  \be\label{cr84}
\frac{\pa u^i}{\pa x^k}=A^i_k+\frac{\pa u^1}{\pa x^1}\delta^i_k,
  \ee
where the matrix $(A^i_k)$ is skew-symmetric. Now by (\ref{cr83})
and (\ref{cr84}) we easily get $\widetilde{r}_{ij}=0$. If $n=2$,
using (\ref{cr81}) we can easily show that $\widetilde{\beta}$ is
closed. Then plus $\widetilde{r}_{ij}=0$, $\widetilde{\beta}$ is
parallel with respect to $\widetilde{\alpha}$.  \qed

\bigskip

It is shown in \cite{Y} that if $n\ge 3$, then all the solutions
to (\ref{cr81}) are given by
 \be\label{cr85}
 u^i=-2\big(\lambda+\langle
 e,x\rangle\big)x^i+|x|^2e^i+q_k^ix^k+f^i,
 \ee
where $\lambda$ is a constant number, $e,f$ are constant
$n$-vectors and
 the constant matrix $(q_k^i)$ is skew-symmetric.
For simplicity, let $(q_k^i)=0$ and $e=tf$ for some constant $t$
in (\ref{cr85}). It is easy to verify that $\widetilde{\beta}$
determined by (\ref{cr80}) and (\ref{cr85}) is closed.  Then by
Lemma \ref{lem71},  $\widetilde{\beta}$ is parallel with respect
to $\widetilde{\alpha}$. Further, we can verify that if $tf\ne 0$,
then $\widetilde{\alpha}$ is of constant sectional curvature if
and only if $\lambda^2+t|f|^2=0$. In this case,
$\widetilde{\alpha}$ is flat.

\begin{ex}
 Defined $\alpha$ and $\beta$ by (\ref{ycw17}), where $\widetilde{\alpha}$ and
 $\widetilde{\beta}$ are determined by (\ref{cr80}). Let $u$ have the
 following form
  $$
 u^i=-2\big(\lambda+t\langle
 f,x\rangle\big)x^i+t|x|^2f^i+f^i,
  $$
  where $t$ is a constant and $f$ is a constant vector satisfying $tf\ne
  0$ and $\lambda^2+t|f|^2\ne 0$. Then the $m$-Kropina metric
  $F=\alpha^m\beta^{1-m}$ is Douglasian but not locally
  projectively flat, where $m\ne 0, 1$.
\end{ex}

\bigskip

{\bf Acknowledgement:}

The  author expresses his sincere thanks to China Scholarship
Council for its funding support.  He did this job  during the
period (June 2012--June 2013) when he as a postdoctoral researcher
visited Indiana University-Purdue University Indianapolis, USA.

\vspace{0.6cm}

\noindent Guojun Yang \\
Department of Mathematics \\
Sichuan University \\
Chengdu 610064, P. R. China \\
{\it  e-mail :} ygjsl2000@yahoo.com.cn


\begin{thebibliography}{999}
\bibitem{AIM} P.  Antonelli, R. Ingarden and M. Matsumoto,
{\it The theory of sprays and Finsler spaces with applications in
physics and biology}, Kluwer Academic Publishers, 1993.

\bibitem{AHM} P. Antonelli, B. Han and J. Modayil, {\it New results on
2-dimensional constant sprays with an application to
heterochrony}, World Scientific Press, 1998.

\bibitem{BaMa} S. Bacso and M. Matsumoto, {\it On Finsler Spaces
of Douglas type. A generalization of the notion of Berwald space},
Publ. Math. Debrecen, {\bf 51}(1997), 385-406.

\bibitem{Dou} J. Douglas, {\it The general geometry of paths},
Ann. of Math. {\bf 29}(1927-28), 143-168.


\bibitem{LSS} B. Li, Y. B. Shen and Z. Shen, {\it On a class of Douglas metrics in
Finsler geometry}, Studia Sci. Math. Hungarica, {\bf 46}(3),
(2009), 355-365.


\bibitem{LS1} B. Li  and Z. Shen, {\it On a class of
projectively fat Finsler metrics with constant flag curvature},
Int. J. of Math.,  {\bf 18}(7) (2007), 1-12.

\bibitem{LS2} B. Li  and Z. Shen, {\it On a class of
weakly Landsberg metrics}, Science in China, Series A, {\bf
50}(2007), 75-85.



\bibitem{RR} M. Rafie-Rad,  {\it Time-optimal solutions of parallel navigation and
Finsler geodesics}, Nonlinear Anal, RWA, 11(2010), 3809C3814.

\bibitem{Shen1}   Z. Shen, {\it On projectively flat
($\alpha,\beta$)-metrics},  Canadian Math. Bull., {\bf
52}(1)(2009), 132-144.

\bibitem{Shen2}  Z. Shen, {\it  On a class of Landsberg metrics in Finsler geometry},
 Canadian J. of Math., {\bf 61}(6) (2009), 1357-1374.

 \bibitem{Shen3} Z.  Shen, {\it Differential geometry of spray and
 Finsler spaces}, Kluwer Academic Publishers, Dordrecht, 2001.

\bibitem{ShY}  Z. Shen and G. Yang, {\it On a class of weakly
Einstein Finsler metrics}, Israel. J. of Math., 2013. (To appear)

 \bibitem{YN} T. Yajima and H. Nagahama, {\it Zermelos condition and seismic ray path}, Nonlinear Anal. RWA, 8(2007), 130C135.



\bibitem{Y} G. Yang, {\it On Randers metrics of isotropic
S-curvature}, Acta Math. Sin., {\bf 52}(6)(2009), 1147-1156 (in
Chinese).


\bibitem{Y2} G. Yang, {\it On Randers metrics of isotropic
S-curvature II}, Publ. Math. Debreceen, {\bf 78}(1) (2011), 71-87.

\bibitem{Y1} G. Yang, {\it On a class of two-dimensional
Douglas and projectively flat Finsler metrics}, preprint.

\bibitem{Y3} G. Yang, {\it On a class of two-dimensional
singular Douglas and projectively flat Finsler metrics}, preprint.


\bibitem{Y4} G. Yang, {\it On a Class of Singular Projectively Flat Finsler Metrics with
Constant Flag Curvature}, preprint.

\bibitem{Y5} G. Yang, {\it On $m$-Kropina metrics of scalar Flag Curvature}, preprint.


\bibitem{YO} R. Yoshikawa and K. Okubo, {\it Kropina spaces of constant curvature II}, arXiv: math/1110.5128v1 [math.DG]
24 Oct 2011.

\bibitem{Yu} C. Yu, {\it A new class of metric deformations and their applications in Finsler geometry}, Ph.D. thesis, Beijing University, (2009).

\end{thebibliography}
\end{document}